\documentclass[11pt]{article}
\usepackage{latexsym,amssymb,amsmath,amsfonts,amsthm}
\makeatletter
\newif\ifmsbmloaded@
\def\loadmsbm{\msbmloaded@true
  \font\tenmsb=msbm10 scaled 1\@ptsize00
  \font\sevenmsb=msbm7 scaled 1\@ptsize00
  \font\fivemsb=msbm5 scaled 1\@ptsize00
  \alloc@8\fam\chardef\sixt@@n\msbfam
  \textfont\msbfam=\tenmsb
  \scriptfont\msbfam=\sevenmsb
  \scriptscriptfont\msbfam=\fivemsb
  }
\def\nonmatherr@#1{\errmessage%
{LateX error: \string#1\space allowed only in math mode}}
\def\Bbb{\relax\ifmmode\expandafter\Bbb@\else
  \expandafter\nonmatherr@\expandafter\Bbb\fi}
\def\Bbb@#1{{\Bbb@@{#1}}}
\def\Bbb@@#1{\fam\msbfam\relax#1}
\@addtoreset{equation}{section}

\makeatother \loadmsbm
\def\R{\Bbb R}

\def\N{\Bbb N}
\def\Z{\Bbb Z}

\def\pa{\partial}

\def\na{\nabla}

\def\sqr#1#2{\vbox{\hrule height .#2pt
\hbox{\vrule width .#2pt height #1pt \kern #1pt \vrule width
.#2pt}\hrule height .#2pt }}
\def\square{\sqr74}
\def\endproof{\hphantom{MM}\hfill\llap{$\square$}\goodbreak}

\newcommand{\beq}{\begin{equation}}
\newcommand{\eeq}{\end{equation}}
\newcommand{\ben}{\begin{eqnarray}}
\newcommand{\een}{\end{eqnarray}}
\newcommand{\beno}{\begin{eqnarray*}}
\newcommand{\eeno}{\end{eqnarray*}}
\newtheorem{Theorem}{Theorem}[section]
\newtheorem{Lemma}[Theorem]{Lemma}

\newtheorem{Remark}{Remark}[section]

\begin{document}
\title{The Beale-Kato-Majda  criterion to the 3D Magneto-hydrodynamics equations }
\date{}
\author{Qionglei Chen $^1$,
Changxing Miao $^1$,
Zhifei Zhang $^2$\\
{\small  $^{1}$Institute of Applied Physics and Computational Mathematics,}\\
    {\small P.O. Box 8009, Beijing 100088, P. R. China.}\\
     {\small (chen\_qionglei@iapcm.ac.cn  and   miao\_changxing@iapcm.ac.cn)}\\
{\small $^{2}$ School of Mathematical Science, Peking University,}\\
{\small Beijing 100871, P. R.  China.}\\
{\small(zfzhang@math.pku.edu.cn)  }}

\maketitle
 \vspace{-1.2in} \vspace{.9in} \vspace{0.2cm}

{\bf Abstract.} We study the blow-up criterion of
smooth solutions to the 3D MHD equations. By means of
the Littlewood-Paley
decomposition, we prove a Beale-Kato-Majda type
 blow-up criterion of smooth
solutions via the vorticity of velocity only, i. e. $\sup_{j\in\Z}\int_0^T\|\Delta_j(\na\times u)\|_\infty dt$, where
$\Delta_j$ is a frequency localization on $|\xi|\approx 2^j$.
 \vspace{0.2cm}

 \vspace{0.2cm}


\section{Introduction}
\label{intro} We consider the 3D incompressible magneto-hydrodynamics (MHD) equations
\begin{align} \label{1.1}(\rm MHD)\,\, \left\{
\begin{aligned}
&\frac{\pa u}{\pa t}-\nu\Delta u+u\cdot \nabla u=-\nabla p-\frac{1}{2}\nabla b^2+b\cdot \nabla
b,\\
&\frac{\pa b}{\pa t}-\eta\Delta b+u\cdot \nabla b=b\cdot \nabla u ,\\
&\nabla\cdot u=\nabla\cdot b=0,\\
&u(0,x)=u_0(x),\quad b(0,x)=b_0(x).
\end{aligned}
\right. \end{align}
Here $u$, $b$ describe
the flow velocity vector and the magnetic field vector
respectively, $p$ is a scalar pressure, $\nu>0$ is the kinematic viscosity and $\eta>0$ is the magnetic diffusivity, while $u_0$ and $b_0$ are
 the given initial velocity and initial magnetic field respectively,
with $\nabla\cdot u_0=\nabla\cdot b_0=0$. If $\nu=\eta=0$, (\ref{1.1}) is called the ideal MHD equations.

Using the standard energy method, it can be easily proved that for given initial data $(u_0, b_0)\in H^s(\R^3)$ with $s>\frac 1 2$,
there exists a positive time $T=T(\|(u_0,b_0)\|_{H^s})$ and a unique smooth solution $(u(t,x),b(t,x))$ on $[0,T)$ to the
MHD equations  satisfying
$$(u, b)\in C([0,T);H^s)\cap C^1((0,T);H^{s})\cap C((0,T); H^{s+2}).$$
Whether  smooth solutions of (\ref{1.1}) on $[0,T)$ will lead to a singularity at $t=T$ is an outstanding open problem, see Sermange and Temam \cite{Ser}.
Caflisch, Klapper and Steele \cite{Caf}
extended the well-known result  of Beale-Kato-Majda \cite{Bea} for the incompressible  Euler equations
to the 3D ideal MHD equations, precisely, under the condition:
\begin{align}\label{1.2}
\int_{0}^{T}(\|\nabla\times u(t)\|_{\infty}+\|\nabla\times b(t)\|_{\infty})dt<\infty
,
\end{align}
then smooth solutions $(u,b)$ can be extended beyond $t=T$.
Recently, there are some researches which have refined (\ref{1.2}) such as
\begin{align}
&\int_{0}^{T}(\|\nabla\times u(t)\|_{\dot{B}^0_{\infty,\infty}}+\|\nabla\times b(t)\|_{\dot{B}^0_{\infty,\infty}}) dt<\infty
,\qquad (\mbox{see}\,\cite{ZL}).\nonumber\\
&\lim_{\varepsilon\rightarrow0}\sup_{j\in\Z}\int_{T-\varepsilon}^{T}
(\|\Delta_j(\nabla\times u)(t)\|_{\infty}+\|\Delta_j(\nabla\times b)(t)\|_{\infty}) dt=\delta<M
,\qquad (\mbox{see}\,\cite{CCM})\nonumber
\end{align}
for some positive constant $M$, and $\Delta_j$ is a frequency
localization on $|\xi|\approx 2^j$. These results can be easily
extended to (\ref{1.1}) with $\nu, \eta>0$. Wu \cite{Wu3} also
extended some Serrin type criterion for the Navier-Stokes equations
to the MHD equations. Many relevant results can be found in
\cite{Wu1,Wu2} and references therein.

However, some numerical experiments \cite{Ha,PP} seem to indicate
that the velocity field plays the more important role than the
magnetic field in the regularity theory of solutions to the MHD
equations. Recently, inspired by Constantin and Fefferman initial work \cite{Con} where
the regularity condition of the direction of vorticity was used to describe
the regularity criterion to the Navier-Stokes equations,  He and Xin \cite{HX} extended it to the  MHD equations,
but did not impose any condition on the magnetic field $b$ which was consistent with the result of numerical experiments.
Precisely, they
showed that the solution remains smooth on $[0,T]$ if the vorticity of the velocity
$w=\na\times u$ satisfies  the following condition
\begin{align}\label{1.3}
\big|w(x+y, t)-w(x, t)\big|\le K|w(x+y,t)||y|^{\frac12}\quad \mbox{if} \,\,\,|y|\le\rho\quad|w(x+y,t)|\ge\Omega,
\end{align}
for $t\in[0,T]$ and three positive constants $K$, $\rho$, $\Omega$.
Also, they \cite{HX} and Zhou \cite{Zhou} obtained some integrability condition of the magnitude of  the only velocity $u$ alone, or  the only
gradient of the velocity $\na u$ alone to characterize the
 regularity criterion to the MHD equations, i.e.
\begin{align}\label{1.4}
\int_0^T \|u(t)\|^q_{p}dt<\infty,\quad\,  \frac{2}{q}+\frac{3}{p}\le1\quad 3<p\le\infty;
\end{align}
or
\begin{align}\label{1.5}
\int_0^T\|\na u(t)\|_{p}^qdt<\infty,\quad\,  \frac{2}{q}+\frac{3}{p}\le2 \quad \frac32<p\le\infty.
\end{align}
We restrict ourselves to substitute $\na u$ by the vorticity $w$ in (\ref{1.5}).
In the case $p<\infty$, using  the Biot-Savart law \cite{Maj1} and the bounds of the Riesz transforms \cite{Ste}
on $L^p$$(1<p<\infty)$, the condition (\ref{1.5}) can be replaced by
\begin{align}\label{1.6}
\int_0^T\|\na \times u(t)\|_{p}^qdt<\infty,\quad \frac{2}{q}+\frac{3}{p}\le2, \quad \frac32<p<\infty.
\end{align}
However, since the lack of continuity of Riesz transforms on $L^\infty$,
their results missed the important marginal case $p=\infty$ which exactly corresponds to the  Beale-Kato-Majda  criterion.
In the case of the Euler equations, Beale, Kato, and Majda get around this difficulty by using the logarithmic Sobolev inequality:
\beq
\|\na u\|_\infty\le C(1+\|\nabla\times u\|_\infty\log(e+\|u\|_{H^s})), \quad s>5/2.\label{1.7}
\eeq
For a refined form of this inequality, it can be referred to \cite{Koz2,Oga}.
In order to make use of  (1.7), one need to estimate the higher order derivatives of the solution
(at least in $H^s, s>5/2$). But, in the case of the MHD equations, it seems difficult to control them
by  the only $\|\na \times u\|_\infty$. Therefore, as in \cite{Caf}, if the logarithmic Sobolev inequality (\ref{1.7}) is used, one can only
derive a criterion described by the vorticity of $u$ and $b$  .
This difficulty is avoided by the following two observations.
On one hand, the $H^1$ norm of the solution can be used to control any $H^s$ norm of the solution,
while the $H^1$ norm of the solution $(u, b)$ can be controlled by $\|\na u\|_\infty$. On the other hand,
if we make use of Littlewood-Paley decomposition to decompose the nonlinear terms into three parts:
low frequency, middle frequency and high frequency, and deal with  each  part  by virtue of different estimates,
we can refine $\|\na u\|_\infty$
to $\|\na \times u\|_\infty$. It should be pointed out that we do not apply the Littlewood-Paley decomposition
to the equation itself as some researches do before, since if we localize the equation on a dyadic partition,
we cannot control the $H^1$ norm of the solution $(u, b)$ via $\|\na\times u\|_{L^\infty}$ in the end when summing up every dyadic partition.

Finally, we remark that the blow-up criterion we will establish  in the framework of
mixed time-space Besov spaces may be the most relaxed in some sense as for the incompressible Euler equations \cite{Pl}
and  the Idea MHD equations \cite{CCM},
where the losing estimate  for the solution
and the logarithmic Sobolev inequality are applied to set up the blow-up criterion, but in this paper,
if we follow their method, as we mentioned above, we cannot characterize the blow-up of smooth solutions by  $\na\times u$  only.

Now we state our result as follows.

\begin{Theorem}
Let $(u_0,b_0)\in H^s$, $s> \frac{1}{2}$ with $\na\cdot u_0=\na\cdot b_0=0$. Suppose that $(u,b)\in C([0,T);H^s)\cap
C^1((0,T);H^{s})\cap C((0,T); H^{s+2})$ is the smooth solution to (\ref{1.1}).
If there exists an absolute constant $M>0$ such that if
\begin{align}\label{1.8}
\lim_{\varepsilon\rightarrow0}\sup_{j\in\Z}\int_{T-\varepsilon}^T
\|\Delta_j(\nabla\times u)\|_{\infty}dt=\delta<M,
\end{align}
then $\delta=0$, and the solution $(u,b)$ can be extended past time
$t=T$. In other words, if
\begin{align}\label{1.9}
\lim_{\varepsilon\rightarrow0}\sup_{j\in\Z}\int_{T-\varepsilon}^T
\|\Delta_j(\nabla\times u)\|_{\infty}dt\ge M,
\end{align}
then the solution blows up at $t=T$. Here $\Delta_j$ is a frequency localization on $|\xi|\approx 2^j$, see Section 2.
\end{Theorem}

\begin{Remark} For the Navier-Stokes equations (with $b=0$ in (\ref{1.1})), Kozono, Taniuchi \cite{Koz1}, and Kozono, Ogawa, Taniuchi \cite{Koz2}
refined the Beale-Kato-Majda criterion to
$$
\int_0^T\|\nabla\times u(t)\|_{BMO}dt<\infty\quad \hbox{and}\quad \int_0^T\|\nabla\times u(t)\|_{\dot B^0_{\infty, \infty}}dt<\infty,
$$
respectively. Here $\dot B^0_{\infty, \infty}$ stands for the homogenous Besov spaces, see Section 2 for the definition.
The condition (\ref{1.8}) is weaker than
all the above mentioned conditions. Hence, this also improves the results of \cite{Koz1, Koz2}.
For the further explanation of (\ref{1.8}), it can be referred to the remarks after Theorem 1 in \cite{Pl}.
\end{Remark}

\begin{Remark}
Very recently, Wu \cite{Wu4} use  the energy estimate combined with the Bony paraproduct technique to derive many
interesting regularity criterions for the generalized MHD equations. In the marginal case,
the regularity criterion obtained there can be expressed as
$$
\int_0^T\|u(t)\|_{B^1_{\infty,\infty}}^{1+\delta}dt<\infty,\quad or \quad \int_0^T\|u(t)\|_{B^{1+\epsilon}_{\infty,\infty}}dt<\infty,
$$
for some $\delta,\epsilon>0$. Here $B^s_{p,q}$ stands for the inhomogenous Besov spaces, see Section 2 for the definition.
But the velocity $u$ cannot be replaced by its vorticity, since the Riesz transformation is not bounded in $B^s_{\infty,\infty}$.
\end{Remark}

\begin{Remark}
For the Ideal MHD equations, whether similar result holds is still open, since the viscous term plays an important role in our proof.
\end{Remark}

\noindent{\bf Notation:} Throughout the paper, $C$ stands for a
``harmless" constant, and changes from line to line; $\|\cdot\|_{p}$ denotes the norm of the Lebesgue space $L^p$.

\section{Preliminaries}
Let ${\cal S}(\R^3)$ be
the Schwartz class of rapidly decreasing functions. Given $f\in
{\cal S}(\R^3)$, its Fourier transform ${\cal F}f=\hat f$ is
defined by
$$
\hat f(\xi)=(2\pi)^{-\frac{3}2}\int_{\R^3}e^{-ix\cdot \xi}f(x)dx.
$$
Now let us recall the Littlewood-Paley
decomposition (see  \cite{Ch1, Tri}).
Choose two nonnegative radial functions $\chi$, $\varphi \in {\cal
S}(\R^3)$, supported respectively in ${\cal B}=\{\xi\in\R^3,\,
|\xi|\le\frac{4}{3}\}$ and ${\cal C}=\{\xi\in\R^3,\,
\frac{3}{4}\le|\xi|\le\frac{8}{3}\}$ such that
\beno
\chi(\xi)+\sum_{j\ge0}\varphi(2^{-j}\xi)=1,\quad\xi\in\R^3,\\
\sum_{j\in\Z}\varphi(2^{-j}\xi)=1,\quad\xi\in\R^3\backslash \{0\}.
\eeno
Let $h={\cal F}^{-1}\varphi$ and $\tilde{h}={\cal F}^{-1}\chi$,
the frequency localization operator is defined by
\beno
&&\Delta_jf=\varphi(2^{-j}D)f=2^{3j}\int_{\R^3}h(2^jy)f(x-y)dy, \\
&&S_jf=\chi(2^{-j}D)f=2^{3j}\int_{\R^3}\tilde{h}(2^jy)f(x-y)dy.
\eeno
Informally, $\Delta_j$ is a
frequency projection to the annulus $\{|\xi|\thickapprox 2^j\}$, while
$S_j$ is a frequency projection to the ball $\{|\xi|\lesssim
2^j\}$. Observe that $\Delta_j=S_j-S_{j-1}$. Also, if $f$ is an
$L^2$ function then $S_jf\rightarrow 0$ in $L^2$ as
$j\rightarrow-\infty$ and $S_jf\rightarrow f$ in $L^2$ as
$j\rightarrow+\infty$(this is an easy consequence of Parseval's
theorem). By telescoping the series, we thus have the homogeneous
Littlewood-Paley decomposition
\begin{equation}\label{2.1}f=\sum_{j=-\infty}^{+\infty}\Delta_jf,\end{equation}
for all $f\in L^2$, where the summation is in the $L^2$ sense.

Let
$s\in \R, 1\le p,q\le\infty$, the homogenous Besov space $\dot
{B}^s_{p,q}$ is defined by
$$\dot {B}^s_{p,q}=\{f\in {\cal Z}'(\R^3); \|f\|_{\dot
{B}^s_{p,q}}<\infty\},$$ where
$$\|f\|_{\dot{B}^s_{p,q}}=
\displaystyle\bigg(\sum_{j=-\infty}^{\infty}2^{jsq}\|\Delta_j f\|_p^q\bigg)^{\frac 1
q},
$$(usual modification if $q=\infty$),
and ${\cal Z}'(\R^3)$  can be identified by the quotient
space of ${\cal S}'/{\cal P}$ with the polynomials space ${\cal
P}$. The inhomogenous Besov space $
{B}^s_{p,q}$ is defined by
$${B}^s_{p,q}=\{f\in {\cal S}'(\R^3); \|f\|_{
{B}^s_{p,q}}=\|S_0f\|_p+\|\{2^{js}\|\Delta_j f\|_p\}_{j\ge 0}\|_{\ell^q}<\infty\}.$$
We now denote  the operator $(I-\Delta)^{\frac{1}{2}}$ by $\Lambda$ which is defined by
$$\widehat{\Lambda f}(\xi)=(1+|\xi|^2)^{\frac{1}{2}}\hat{f}(\xi).$$
More generally, $\Lambda^s f$ for $s\in\R$ can be identified with the Fourier Transform
$$\widehat{\Lambda^s f}(\xi)=(1+|\xi|^2)^{\frac{s}{2}}\hat{f}(\xi).$$
For $s\in\R$, we define
$$\|f\|_{{H}^s}\triangleq\|\Lambda^s f\|_{L^2}\triangleq\bigg(\int_{\R^3}(1+|\xi|^s)^{2}|\hat{f}(\xi)|^2d\xi\bigg)^{\frac 12},$$
and the Sobolev space $H^s$ is denoted  by ${H}^s\triangleq\{f\in {\cal S}'(\R^3); \|f\|_{
{H}^s}<\infty\}.$ The usual Sobolev space $H^{s,p}$ is endowed with the norm
$$\|f\|_{H^{s,p}}\triangleq\|\Lambda^s f\|_{L^p}.$$
We can refer to \cite{Tri} for more details.

\begin{Lemma}
Let $k\in \N$. There exist
constants $C$ independent of $f$, $j$ such that for all  $1\le p\le q\le\infty$
\begin{align}
&\sup_{|\alpha|=k}\|\partial^\alpha\Delta_j f\|_q\le C2^{jk+3j(\frac{1}{p}-\frac{1}{q})}\|\Delta_jf\|_{p},
\label{2.2}\\
&\|R_k\Delta_jf\|_{q}\le C2^{3j(\frac{1}{p}-\frac{1}{q})}\|\Delta_jf\|_p.\label{2.3}
\end{align}
Here $R_k$ $(k=1,2,3)$ is the Riesz transform in $\R^3$.
\end{Lemma}
The proof of this lemma can be found in \cite{Ch1,Mey}.

\begin{Remark}
Suppose that the vector function $f$ is divergence-free, and set $g=\nabla\times f$. Then there exist
constants $C$ independent of $f$ such that
\begin{align}\label{2.4}
\|\na f\|_p\le C\|g\|_p,\quad\forall\,\, 1<p<\infty.
\end{align}
If the frequency of $f$ is restricted to some annulus $\{|\xi|\approx 2^j\}$, then there holds
\begin{align}\label{2.5}
\|\na f\|_p\le C\|g\|_p,\quad\forall\,\, 1\le p\le\infty.
\end{align}
Indeed, the inequality (\ref{2.4}) can be derived from the Biot-Savart law \cite{Maj1}
and the bounds of the Riesz transforms \cite{Ste} on $L^p$$(1<p<\infty)$, while  the inequality (\ref{2.5})
can be deduced from  the Biot-Savart law and (\ref{2.3}) .
\end{Remark}

\begin{Lemma}[Commutator estimate]
Let $1<p<\infty$, $s>0$. Assume that $f, g\in H^{s,p}$,
then there exist
constants $C$ independent of $f$, $g$ such that
\begin{align}\label{2.6}
\|\Lambda^s(fg)-f\Lambda^s g\|_{L^p}\le C(\|\na f\|_{L^{p_1}}\|g\|_{{H}^{s-1,p_2}}+\|f\|_{{H}^{s,p_3}}\|g\|_{L^{p_4}})
\end{align}
with $p_2, p_3\in(1,+\infty)$ such that
$$\frac{1}{p}=\frac{1}{p_1}+\frac{1}{p_2}=\frac{1}{p_3}+\frac{1}{p_4}.$$
\end{Lemma}
This lemma is well-known and for a proof, see \cite{KP}.

\section{Proof of Theorem 1.1}

We will divide the proof of Theorem 1.1 into two steps.\vspace{0.2cm}

{\bf Step 1.}\, $H^1$ estimates. \vspace{0.1cm}

In this step we will show  there exists  $\varepsilon>0$,
\begin{align}\label{3.1}
\sup_{t\in[T-\varepsilon, T)}\big(\|u(t)\|_{H^1}+\|b(t)\|_{H^1}\big)\le C
\big(\|u(T-\varepsilon)\|_{H^1}+\|b(T-\varepsilon)\|_{H^1}+e\big).
\end{align}
Let $w(t,x)=\na\times u(t,x)$ and $J(t,x)=\na \times b(t,x)$. Taking the curl on both sides of (\ref{1.1}),
it can  be written as
\begin{align}\label{3.2}\left\{
\begin{aligned}
&\frac{\pa w}{\pa t}-\nu\Delta w+(u\cdot\na)w-(w\cdot\na)u-(b\cdot\na)J+(J\cdot\na)b=0,\\
&\frac{\pa J}{\pa t}-\eta\Delta J+(u\cdot\na)J-(J\cdot\na)u-(b\cdot\na)w+(w\cdot\na) b=2T(b, u)
\end{aligned}\right.\end{align}
with $$T(b,u)=\left(\begin{array}{ll}
\pa_2b\cdot\pa_3 u-\pa_3b\cdot\pa_2u\\
\pa_3b\cdot\pa_1 u-\pa_1b\cdot\pa_3 u\\
\pa_1b\cdot\pa_2 u-\pa_2b\cdot\pa_1 u
\end{array}\right).$$
Multiplying the first equation of (\ref{3.2}) by $w$, the second one of (\ref{3.2}) by $J$,  then
adding the resulting equations yields that
\begin{align}\label{3.3}
&\frac{1}{2}\frac{d}{dt}(\|w(t)\|_2^2+\|J(t)\|_2^2)+\nu\|\na
w(t)\|_2^2+\eta\|\na J(t)\|_2^2=\int_{\R^3}(w\cdot\na) u\cdot wdx+
\int_{\R^3}(J\cdot\na) u\cdot Jdx\nonumber\\
&\qquad-\int_{\R^3}\big((J\cdot\na) b\cdot w+ (w\cdot\na) b\cdot J\big)dx+2\int_{\R^3}T(b,u)\cdot Jdx\nonumber\\
&\triangleq I+II+III+2\,IV,
\end{align}
where we have used the facts $$\int_{\R^3}(u\cdot\na) w\cdot wdx=\int_{\R^3}(u\cdot\na) J\cdot Jdx=0$$ and
$$\int_{\R^3}\big( (b\cdot\na) J\cdot w+(b\cdot\na) w\cdot J\big)dx=0,$$
which can be deduced from the  ${\rm div}u={\rm div}b=0$ and integrating by part.

In what follows, we will deal with each term on the right hand side of (\ref{3.3})
separately below. Let us begin with estimating the term $I$. Using the Littlewood-Paley decomposition (\ref{2.1}), we have
\begin{align}\label{3.4}
w=\sum_{j\in\Z}\Delta_j w=\sum_{j<-N}\Delta_jw+\sum_{-N\le j\le N}\Delta_j w+\sum_{j>N}\Delta_jw,
\end{align}
where $N$ is a positive integer to be determined later. Putting (\ref{3.4}) into $I$ produces that
\begin{align}
&I=\sum_{j<-N}\int_{\R^3}(w\cdot\na) u\cdot\Delta_jw dx+\sum_{j=-N}^N\int_{\R^3}(w\cdot\na) u\cdot\Delta_jw dx
+\sum_{j>N}\int_{\R^3} (w\cdot\na) u\cdot\Delta_jw dx\nonumber\\&\triangleq I_1+I_2+I_3.\nonumber
\end{align}
Using the H\"{o}lder inequality, (\ref{2.4}) and (\ref{2.2}) to obtain that
\begin{align}\label{3.5}
&|I_1|\le\|w\|_2\|\na u\|_2\sum_{j<-N}\|\Delta_jw\|_{\infty}\le C\|w\|_2^2\sum_{j<-N}2^{\frac{3}{2}j}\|\Delta_j w\|_2\le C2^{-\frac{3}{2}N}\|w\|_2^3,\\
&|I_2|\le \|w\|_2\|\na u\|_2\sum_{-N\le j\le N}\|\Delta_j w\|_\infty\le C\|w\|_2^2\sum_{-N\le j\le N}\|\Delta_jw\|_\infty.
\end{align}
From the H\"{o}lder inequality, (\ref{2.4}), (\ref{2.2}) and Gagliardo-Nirenberg inequality, it follows that
\begin{align}\label{3.7}
|I_3|&\le\|w\|_3\|\na u\|_3\sum_{j>N}\|\Delta_jw\|_3\le C\|w\|_3^2\sum_{j>N}2^{\frac j 2}\|\Delta_jw\|_2\nonumber\\&
\le C\|w\|_3^2\bigg(\sum_{j>N}2^{-j}\bigg)^{\frac 1 2}\bigg(\sum_{j>N}2^{2j}\|\Delta_jw\|_2^2\bigg)^{\frac 1 2}\le C2^{-\frac N 2}\|w\|_2\|\na w\|_2^2.
\end{align}
By summing up (\ref{3.5})-(\ref{3.7}), we get
\begin{align}\label{3.8}
|I|\le C\big(2^{-\frac{3}{2}N}\|w\|_2^3+\|w\|_2^2\sum_{-N\le j\le N}\|\Delta_jw\|_\infty+2^{-\frac N 2}\|w\|_2\|\na w\|_2^2\big).
\end{align}
Using the Littlewood-Paley decomposition (\ref{2.1}) to $\na u$,  $II$ can be written as
\begin{align}
II=\sum_{j<-N}\int_{\R^3} (J\cdot\na) \Delta_ju\cdot Jdx+\sum_{j=-N}^{N}\int_{\R^3} (J\cdot\na) \Delta_ju\cdot Jdx
+\sum_{j>N}\int_{\R^3} (J\cdot\na) \Delta_ju\cdot Jdx.\nonumber
\end{align}
Then the H\"{o}lder inequality, (\ref{2.2}),  (\ref{2.5}) and Gagliardo-Nirenberg inequality allow us to show that
\begin{align}\label{3.9}
|II|\le C\big(2^{-\frac{3}{2}N}\|w\|_2\|J\|_2^2+\|J\|_2^2\sum_{-N\le j\le N}\|\Delta_jw\|_\infty+2^{-\frac N 2}\|J\|_2\|\na w\|_2\|\na J\|_2\big).
\end{align}
Similarly, using the Littlewood-Paley decomposition (\ref{2.1}), $III$ and $IV$ can be written respectively as
\begin{align}
&III=\sum_{j<-N}\int_{\R^3}J\cdot\big(\na b+(\na b)^T\big)\cdot \Delta_jwdx+\sum_{j=-N}^{N}\int_{\R^3}J\cdot \big(\na b+(\na b)^T\big)\cdot\Delta_jwdx
\nonumber\\&\qquad+\sum_{j>N}\int_{\R^3}J\cdot \big(\na b+(\na b)^T\big)\cdot\Delta_jwdx,\nonumber\\
&IV=\sum_{j<-N}\int_{\R^3}T(b, \Delta_ju)\cdot Jdx+\sum_{j=-N}^{N}\int_{\R^3}T(b, \Delta_ju)\cdot Jdx
+\sum_{j>N}\int_{\R^3}T(b, \Delta_ju)\cdot Jdx.
\nonumber
\end{align}
Then exactly as in the derivation of (\ref{3.8}), (\ref{3.9}), we can deduce that
\begin{eqnarray}\label{3.10}
|III|+2|IV|&\le& C\big(2^{-\frac{3}{2}N}\|w\|_2\|J\|_2^2+\|J\|_2^2\sum_{-N\le j\le N}\|\Delta_jw\|_\infty\nonumber\\
&&+2^{-\frac N 2}\|J\|_2\|\na w\|_2\|\na J\|_2\big).
\end{eqnarray}
Combining (\ref{3.8})-(\ref{3.10}) with (\ref{3.3}), Young inequality yields that for $t\in[0,T)$
\begin{align}\label{3.11}
&\frac{d}{dt}\big(\|w(t)\|_2^2+\|J(t)\|_2^2\big)+2\nu\|\na w(t)\|_2^2+2\eta\|\na J(t)\|_2^2\nonumber\\&\le C\bigg(
2^{-\frac{3}{2}N}\big(\|w(t)\|_2^3+\|J(t)\|_2^3\big)+\sum_{-N\le j\le N}
\|\Delta_jw(t)\|_\infty\big(\|w(t)\|_2^2+\|J(t)\|_2^2\big)\nonumber\\&\quad+
2^{-\frac N 2}\big(\|w(t)\|_2+\|J(t)\|_2\big)\big(\|\na w(t)\|_2^2+\|\na J(t)\|_2^2\big)\bigg).
\end{align}
Now let us choose a fixed positive integer $N$ such that $C2^{-\frac N 2}\big(\|w(t)\|_2+\|J(t)\|_2\big)\le\min(\nu,\eta)$, i.e.
\begin{align}\label{3.12}
N\ge\bigg[\frac{2}{\log 2}\log^+\bigg(\frac{C}{\min(\nu,\eta)}\big(\|w(t)\|_2+\|J(t)\|_2\big)\bigg)\bigg]+1.
\end{align}
where $\log^+x=\log(e+x)$. Thus (\ref{3.11}) and (\ref{3.12}) imply that for  $t\in[0, T)$
\begin{align}\label{3.13}
&\frac{d}{dt}\big(\|w(t)\|_2^2+\|J(t)\|_2^2\big)+\nu\|\na w(t)\|_2^2+\eta\|\na J(t)\|_2^2\nonumber\\&
\le C\sum_{j=-N}^N\|\Delta_jw(t)\|_\infty\big(\|w(t)\|_2^2+\|J(t)\|_2^2\big)+C,
\end{align}
which together with the Gronwall inequality gives that for $t\in[0, T)$
\begin{align}
\|w(t)\|_2+\|J(t)\|_2\le\exp\bigg(C\sum_{j=-N}^N\int_0^t\|\Delta_jw(t')\|_\infty dt'\bigg)(\sqrt{Ct}+\|w(0)\|_2+\|J(0)\|_2).\nonumber
\end{align}
Recalling the choice of $N$ in (\ref{3.12}), it follows from the above estimate that
\begin{align}
\|w(t)\|_2+\|J(t)\|_2&\le\exp\bigg(C\log^+\big(\|w(t)\|_2+\|J(t)\|_2\big)\sup_{j\in\Z}\int_0^t\|\Delta_jw(t')\|_\infty dt'\bigg)
\nonumber\\&\quad\times (\sqrt{Ct}+\|w(0)\|_2+\|J(0)\|_2).\nonumber
\end{align}
For simplicity, let $E(t)\triangleq\|w(t)\|_2+\|J(t)\|_2,$
the above inequality implies that
\begin{align}\label{3.14}
\sup_{[0,T)}E(t)\le \exp\bigg(C\log^+\big(\sup_{[0,T)}E(t)\big)\sup_{j\in\Z}\int_0^T\|\Delta_jw(t')\|_\infty dt'\bigg)
(\sqrt{CT}+E(0)).
\end{align}
We point out that the inequality (\ref{3.14}) still holds if the time interval is replaced by $[T-\varepsilon, T)$.
It follows from (3.14) that
\begin{align}
\sup_{t\in [T-\varepsilon,T)}E(t)&\le\exp\bigg(\log^+\sup_{t\in [T-\varepsilon,T)}E(t)\sup_{j\in\Z}\int_{T-\varepsilon}^T
\|\Delta_jw(t')\|_\infty dt'\bigg)\nonumber\\&
\quad\times(\sqrt{C\varepsilon}+\|w(T-\varepsilon)\|_2+\|J(T-\varepsilon)\|_2).\nonumber
\end{align}
Defining $Z(T)\triangleq\log\big(\displaystyle\sup_{t\in [T-\varepsilon,T)}E(t)+e\big)$, the above estimate  means that
\begin{align}\label{3.15}
Z(T)\le \log\big(\sqrt{C\varepsilon}+E(T-\varepsilon)+e\big)+CZ(T)\sup_{j\in\Z}\int_{T-\varepsilon}^T\|\Delta_j w(t')\|_\infty dt'.
\end{align}
If we choose $M=\frac{1}{2C}$ in Theorem 1.1, the condition (\ref{1.8}) ensures that there exists a small positive $\varepsilon_0$ such that
$$C\displaystyle\sup_{j\in\Z}\int_{T-\varepsilon}^T\|\Delta_jw(t')\|_\infty dt'\le\frac12, \qquad \forall\, \varepsilon\in (0,\varepsilon_0).$$
Then the inequality (\ref{3.15}) implies that
\begin{align}\label{3.16}
Z(T)\le 2\log\big(E(T-\varepsilon)+e\big), \qquad \forall \,\varepsilon\in (0,\varepsilon_0).
\end{align}
On the other hand,  it is easy to prove that
the solution $(u,b)$ satisfies the energy inequality
\begin{align}
\|u(t)\|_2^2+\|b(t)\|_2^2+2\int_s^t\big(\nu\|\na u(t')\|_2^2+\eta\|\na b(t')\|_2^2\big)dt'\le
\|u(s)\|_2^2+\|b(s)\|_2^2,\nonumber
\end{align}
which together with (\ref{3.16}), (\ref{2.4}) implies (\ref{3.1}).

\vspace{.2cm}

{\bf Step 2.} $H^s$$(s>1)$ estimates. \vspace{0.1cm}

For completeness,
we will show how to deduce $H^s$ estimates from $H^1$ estimates.
Taking the operation $\Lambda^s$ on both sides of (\ref{1.1}),
multiplying $(\Lambda^s u, \Lambda^s
b)$ to the resulting equation, and integrating over $\R^3$, we get
\begin{align}\label{3.17} &\frac{1}{2}
\frac{d}{dt}(\|\Lambda^s u(t)\|_2^2+\|\Lambda^s
b(t)\|_2^2)+\nu\|\na\Lambda^s
u(t)\|_2^2+\eta\|\na\Lambda^s
b(t)\|_2^2\nonumber\\&=-\int_{\R^3}\Lambda^s(u\cdot\na
u)\Lambda^s u dx+ \int_{\R^3}\Lambda^s(b\cdot\na
b)\Lambda^s u dx-\int_{\R^3}\Lambda^s(u\cdot\na
b)\Lambda^s bdx\nonumber\\&\quad+ \int_{\R^3}\Lambda^s(b\cdot\na
u)\Lambda^sbdx.\end{align}
Noting that  ${\rm div}u={\rm div}b=0$  and
integrating by parts, we rewrite (\ref{3.17}) as
\begin{align}\label{3.18} &\frac{1}{2}
\frac{d}{dt}(\|\Lambda^s u(t)\|_2^2+\|\Lambda^s
b(t)\|_2^2)+\nu\|\na\Lambda^s
u(t)\|_2^2+\eta\|\na\Lambda^s
b(t)\|_2^2\nonumber\\&=-\int_{\R^3}(\Lambda^s(u\cdot\na
u)-u\cdot\Lambda^s\na u)\Lambda^s u
dx-\int_{\R^3}(\Lambda^s(u\cdot\na b)-u\cdot\Lambda^s\na
b)\Lambda^s
bdx\nonumber\\&\quad+\int_{\R^3}(\Lambda^s(b\cdot\na
b)-b\cdot\Lambda^s\na b)\Lambda^s u+
(\Lambda^s(b\cdot\na u)-b\cdot\Lambda^s\na
u)\Lambda^s
bdx\nonumber\\&\triangleq
\Pi_1+\Pi_2+\Pi_3.\end{align}
For  the  term $\Pi_1$,
it follows from (\ref{2.6}), H\"{o}lder inequality,  Gagliardo-Nirenberg inequality and Young inequality that
\begin{align}\label{3.19}
|\Pi_1|&\le C(\|\na u\|_2\|\na u\|_{H^{s-1,4}}+\|u\|_{H^{s,4}}\|\na u\|_2)\|u\|_{H^{s,4}}
\nonumber\\&\le C\|\na u\|_{2}\|u\|_{H^s}^{\frac12}\|\na u\|_{H^s}^{\frac32}
\le C\|\na u\|_2^4\|u\|_{H^s}^2+\frac{\nu}{2}\|\na u\|_{H^s}^2.
\end{align}
The other terms can be treated in the same way:
\begin{align}\label{3.20}
&|\Pi_2|+|\Pi_3|\le C(\|\na u\|_2\|\na b\|_{H^{s-1,4}}+\|u\|_{H^{s,4}}\|\na b\|_{2})\|b\|_{H^{s,4}}
\nonumber\\&\qquad+C(\|\na b\|_2\|\na b\|_{H^{s-1,4}}+\|b\|_{H^{s,4}}\|\na b\|_{2})\|u\|_{H^{s,4}}\nonumber\\&\qquad+
C(\|\na b\|_2\|\na u\|_{H^{s-1,4}}+\|b\|_{H^{s,4}}\|\na u\|_{2})\|b\|_{H^{s,4}}\nonumber\\&
\le C(\|\na u\|_2^4+\|\na b\|_2^4)(\|u\|_{H^s}^2+\|b\|_{H^s}^2)+\frac{\nu}{2}\|\na u\|_{H^s}^2+\frac{\eta}{2}\|\na b\|_{H^s}^2.
\end{align}
By summing up (\ref{3.19}) and (\ref{3.20}) with (\ref{3.18}), we get
\begin{align}
&\frac{d}{dt}(\|u(t)\|^2_{H^s}+\|b(t)\|^2_{H^s})+\nu\|\na u(t)\|_{H^s}^2+\eta\|\na b(t)\|_{H^s}^2\nonumber\\&
\qquad\le C(\|\na u(t)\|_2^4+\|\na b(t)\|_2^4)(\|u(t)\|_{H^s}^2+\|b(t)\|_{H^s}^2).\nonumber
\end{align}
Then the Gronwall inequality yields that
\begin{align}
&\|u(t)\|^2_{H^s}+\|b(t)\|^2_{H^s}+\int_{0}^t(\nu\|\na u(t')\|_{H^s}^2+\eta\|\na b(t')\|_{H^s}^2)dt'
\nonumber\\&\le C(\|u(0)\|^2_{H^s}+\|b(0)\|^2_{H^s})
\exp\bigg(t\sup_{t'\in [0,t)}\|(u(t'), b(t'))\|_{H^1}^4\bigg).
\end{align}
Hence we have the $H^s$ regularity for the solution at $t=T$ and the solution
can be continued after $t=T$. This completes the proof of Theorem 1.1.

\vspace{.2cm}

\textbf{Acknowledgements} Q. Chen and C. Miao  were partly
supported by the NSF of China (10571016)  and  The Institute of
Mathematical Sciences, The Chinese University of Hong Kong.  Z.
Zhang is supported by the NSF of China(10601002). The authors wish
to thank Prof.Zhouping Xin for stimulating discussion about this
problem.

\end{document}